\documentclass[11pt]{article}
\usepackage{amsfonts}
\usepackage[dvipsnames]{xcolor}
\usepackage{amsmath,amssymb,comment}
\usepackage{verbatim}
\usepackage{colonequals}
\usepackage[utf8]{inputenc}
\usepackage[T1]{fontenc}
\usepackage{hyperref}
\hypersetup{
  colorlinks=true,
  linkcolor=black,
  citecolor=black,
  urlcolor=blue,
  pdftitle={Chemotaxis-consumption interaction: Solvability and asymptotics in general high-dimensional domains},
  pdfauthor={Lankeit, Winkler},
  pdfkeywords={},
  bookmarksopen=true,
}

\newtheorem{theo}{Theorem}[section]
\newtheorem{lem}[theo]{Lemma}
\newtheorem{cor}[theo]{Corollary}

\newtheorem{defi}[theo]{Definition}

\newcommand{\mysection}[1]{\section{#1} \setcounter{equation}{0}}
    \makeatletter
\def\@fnsymbol#1{\ensuremath{\ifcase#1\or *\or \ddagger\or
   \mathsection\or \mathparagraph\or \|\or **\or \dagger\dagger
   \or \ddagger\ddagger \else\@ctrerr\fi}}
    \makeatother
\newcommand{\proof}{{\sc Proof.} \quad}
\newcommand{\proofc}{{\sc Proof} \ }
\newcommand{\be}{\begin{equation} \label}
\newcommand{\ee}{\end{equation}}
\newcommand{\bea}{\begin{eqnarray}\label}
\newcommand{\eea}{\end{eqnarray}}
\newcommand{\bas}{\begin{eqnarray*}}
\newcommand{\eas}{\end{eqnarray*}}
\newcommand{\bit}{\begin{itemize}}
\newcommand{\eit}{\end{itemize}}
\newcommand{\qed}{\hfill$\Box$ \vskip.2cm}
\newcommand{\nn}{\nonumber}
\newcommand{\R}{\mathbb{R}}

\newcommand{\N}{\mathbb{N}}
\newcommand{\pO}{\partial\Omega}

\newcommand{\eps}{\varepsilon}

\newcommand{\wto}{\rightharpoonup}
\newcommand{\wsto}{\stackrel{\star}{\rightharpoonup}}

\newcommand{\io}{\int_\Omega}

\newcommand{\na}{\nabla}
\newcommand{\Del}{\Delta}
\newcommand{\del}{\delta}

\newcommand{\lam}{\lambda}

\newcommand{\pa}{\partial}
\newcommand{\bom}{\overline{\Omega}}
\newcommand{\Om}{\Omega}

\newcommand{\wh}{\widehat}

\newcommand{\hs}{\hspace*}

\newcommand{\vp}{\varphi}
\newcommand{\lbal}{\left\{ \begin{array}{l}}
\newcommand{\lball}{\left\{ \begin{array}{ll}}
\newcommand{\ear}{\end{array} \right.}

\newcommand{\abs}{\\[5pt]}

\newcommand{\tme}{T_{max,\eps}}

\newcommand{\ueps}{u_\eps}
\newcommand{\veps}{v_\eps}

\newcommand{\heps}{h_\eps}

\newcommand{\geps}{g_\eps}
\newcommand{\yeps}{y_\eps}

\newcommand{\Feps}{{\mathcal{F}}_\eps}

\newcommand{\D}{{\mathcal{D}}}
\newcommand{\f}[2]{\frac{#1}{#2}}
\newcommand{\Tmaxe}{\tme}
\newcommand{\ue}{\ueps}
\newcommand{\ve}{\veps}
\newcommand{\norm}[2][]{\|#2\|_{#1}}
\newcommand{\Lom}[1]{L^{#1}(\Om)}
\newcommand{\Wom}[2]{W^{#1,#2}(\Om)}
\newcommand{\upto}{\nearrow}
\newcommand{\kl}[1]{\left(#1\right)}
\newcommand{\Ombar}{\overline{\Om}}
\newcommand{\set}[1]{\left\{#1\right\}}
\begin{document}
\enlargethispage{10mm}
\title{Chemotaxis-consumption interaction:\\
Solvability and 		
asymptotics
in general high-dimensional domains}
\author
{
Johannes Lankeit\footnote{lankeit@ifam.uni-hannover.de, ORCID: \href{https://orcid.org/0000-0002-2563-7759}{0000-0002-2563-7759}}\\
{\small Leibniz Universität Hannover, Institut für Angewandte Mathematik,}\\
{\small Welcfengarten 1, 30167 Hannover, Germany}
 \and
Michael Winkler\footnote{michael.winkler@math.uni-paderborn.de}\\
{\small Universit\"at Paderborn, Institut f\"ur Mathematik,}\\
{\small 33098 Paderborn, Germany} }
\date{}
\maketitle
\begin{abstract}
\noindent 
The basic chemotaxis-consumption model 
\bas
	\lbal
	u_t = \Delta u - \nabla \cdot(u\nabla v), \\[1mm]
	v_t = \Delta v - uv
	\ear
\eas
is considered in 
general, possibly non-convex 
bounded domains of arbitrary spatial dimension.
Global existence of weak solutions is shown, along with eventual smoothness of solutions and their stabilization in the large time limit.
\end{abstract}
%
%
%
%
%
%
\newpage
\section{Introduction}\label{intro}
One of the fundamental classes of chemotaxis systems is formed by those with consumption of the signal, which at the latest since their appearance in the description of bacterial behaviour in fluid experiments \cite{goldstein2005} have been at the focus of interest; for an overview, we refer to the survey \cite{lan_win_depleting}. In this article, we are concerned with their most 
proptotypical representative given by
\be{0}
	\left\{ \begin{array}{ll}	
	u_t = \Del u - \na \cdot (u\na v),
	\qquad & x\in\Om, \ t>0, \\[1mm]
	v_t = \Del v - uv,
	\qquad & x\in\Om, \ t>0, \\[1mm]
	\frac{\pa u}{\pa\nu}=\frac{\pa v}{\pa\nu}=0,
	\qquad & x\in\pO, \ t>0, \\[1mm]
	u(x,0)=u_0(x), \quad v(x,0)=v_0(x),
	\qquad & x\in\Om,
	\end{array} \right.
\ee 
to be posed in a smoothly bounded domain $\Om\subset\R^n$ with $n\ge 1$.
In improving local to global existence results, the second equation herein appears to be beneficial, as it immediately provides $\Lom\infty$ bounds for the signal, the second solution component $v$, in stark contrast to the situation in the classical Keller--Segel-production system addressing contexts in which the signal is produced by individuals. This simplicity is deceptive, however, as $\Lom\infty$ bounds are insufficient to adequately cope with the cross-diffusive term in the first equation for the bacterial population density $u$ and, in fact, recent results regarding 
some variants of \eqref{0} involving different boundary conditions (\cite{dong}, \cite{ahn_win})	
provide indications that consumptive cross-diffusion mechanisms as incorporated in \eqref{0} are even compatible with blow-up. 

From this point of view, it may not be too surprising that 
the existence theory for \eqref{0} seem yet incomplete. In bounded subdomains of the plane,
unique global classical solvability has been known for a long time (\cite[Thm.~1.1(i)]{win_CPDE}), 
while in higher-dimensional domains, a smallness condition on $\|v_0\|_{L^\infty(\Om)}$ 
has so far been necessary for similar conclusions (\cite{tao11}, \cite{baghaei}, \cite{heihoff}).
Large initial data can be treated in the three-dimensional case in the context of weak solutions (\cite[Thm.~1.1(ii)]{win_CPDE}).
In higher dimensions, even global weak solvability is still unknown; 
constructions available in the literature seem to be restricted to certain generalized solutions, substantially weakening the 
notion of solvability by requiring that, inter alia, only renormalized versions $\xi (u)$ of the respective first components solve a weak form of their corresponding PDE for each $\xi $ with $\xi '\in C_0^\infty([0,\infty))$ (\cite{li_yuxiang_EJDE}).

Restrictions more crucial than those on the spatial dimension, however, appear to have been underlying essential parts
of the large time analysis of (\ref{0}). Indeed, while in two-dimensional domains the classical solutions to (\ref{0})
are known to stabilize toward the semi-trivial steady state at the corresponding mass level (\cite{jiang}),
similar conclusions for the possibly less regular solutions known to exist in higher-dimensional situations could
to date only be drawn under the additional technical restriction that $\Om$ be convex (\cite{taowin_consumption},
\cite{win_TRAN}).

The intention of the present note is to indicate that also in such higher-dimensional cases, any such extra requirement 
in fact is superfluous; indeed, we shall see that in smoothly bounded domains of any dimension and arbitrary 
geometry, within a natural concept of weak solvability some global solutions can always be found, and that
all these solutions actually become smooth and stabilize to spatially homogeneous equilibria.

To highlight the methodological novelty underlying our approach, let us recall that in apparently all previous
cases the existence analysis is based on a priori bounds obtained from a quasi-energy structure formally expressed in the
identity
\be{en0}
	\frac{d}{dt} \bigg\{ \io u\ln u + \frac{1}{2} \io \frac{|\na v|^2}{v} \bigg\}
	+ \io \frac{|\na u|^2}{u} 
	+ \io v|D^2\ln v|^2
	+ \frac{1}{2} \io \frac{u}{v} |\na v|^2
	= \frac{1}{2} \int_{\pO} \frac{1}{v} \frac{\pa |\na v|^2}{\pa\nu}.
\ee
(cf.~also Lemma \ref{lem3} below).
Now with regard to the intended result on solvability, the decisive difference lies in the treatment of the cross-diffusive term: Whereas in \cite{win_CPDE} some $L^1$ compactness properties of the nonlinear flux $u\na v$ were obtained from separate integrability features of $u$ and $\na v$, namely from bounds on $\int \int u^{\f{n+2}n}$ (cf. \eqref{6.1}) and on $\int \int |\na v|^4$ (whose combination via Hölder's inequality as in \cite[p.348]{win_CPDE} ceases to work in higher dimensions), presently, we will here rather make use of {\em a priori} information on
\bas
	\int \int \frac{u}{v} |\na v|^2
\eas
that appears as part of the integrated dissipation rate in (\ref{en0}).

Similarly, the announced progress related to large time behavior rests on a refined energy analysis as well: In \cite{taowin_consumption}, some first decay information on $v$ ensuring convergence as $t\to\infty$ was taken from finiteness of $\int_0^\infty\io |\na v|^4$; indeed, under the assumption that $\Omega$ is convex this is implied by (\ref{en0}) due to the fact that then $\frac{\pa |\na v|^2}{\pa\nu}\le 0$ on $\pO\times (0,\infty)$, meaning that (\ref{en0}) expresses a genuine energy inequality
in this case. 
Presently, we instead use some basic decay information of $\na v$ (see (\ref{gradv})) to control the possibly ill-signed 
boundary integral in (\ref{en0}) {\em in terms of some quantities that are known to decay} in the large time limit
(see Lemma \ref{lem44} and Lemma \ref{lem4}).
This will be seen to imply temporally sublinear growth of $\int_0^T\io |u-\mu|$ and thus, since $\io v_t+\mu\io v = \io (\mu-u)v$, smallness of $\f1T\int_0^T\io v$, which due to monotonicity of $L^p$ norms of $v$ will turn out to be a suitable starting point for the proof of its convergence.

{\bf Main results.} \quad
In order to state our results more precisely, let us first specify the definition of solutions to be referred to.
\begin{defi}\label{dw}
  Let $n\ge 1$ and $\Om\subset\R^n$ be open and bounded, and assume that $u_0\in L^1(\Om)$ and $v_0\in L^\infty(\Om)$
  are nonnegative.
  Then by a {\em global weak solution} of (\ref{0}) we mean a pair $(u,v)$ of nonnegative functions
  \be{w1}
	\lbal
	u\in L^1_{loc}([0,\infty);W^{1,1}(\Om))
	\qquad \mbox{and} \\[1mm]
	v\in L^\infty_{loc}(\bom\times [0,\infty)) \cap L^1_{loc}([0,\infty);W^{1,1}(\Om))
	\ear
  \ee
  such that 
  \be{w2}
	u\na v \in L^1_{loc}(\bom\times [0,\infty);\R^n),
  \ee
  and that for each $\vp\in C_0^\infty(\bom\times [0,\infty))$, the identities
  \be{wu}
	- \int_0^\infty \io u\vp_t - \io u_0\vp(\cdot,0)
	= - \int_0^\infty \io \na u\cdot\na\vp 
	+ \int_0^\infty \io u\na v\cdot\na\vp
  \ee
  and 
  \be{wv}
	\int_0^\infty \io v\vp_t + \io v_0\vp(\cdot,0)
	= \int_0^\infty \io \na v\cdot\na\vp
	+ \int_0^\infty \io uv\vp
  \ee
  hold.
\end{defi}
In this framework, we then indeed obtain the following.
\begin{theo}\label{theo15}
  Let $n\ge 1$ and $\Om\subset\R^n$ be a bounded domain with smooth boundary,
  and suppose that
  \be{init}
	\lbal
	u_0 \in L\log L(\Om)\setminus\set0 \mbox{ is nonnegative in $\Om$, \quad and that} \\[1mm]
	v_0>0 \text{ is such that } \sqrt{v_0}\in \Wom12\cap\Lom\infty.
	\ear
  \ee
  Then one can find nonnegative functions
  \be{reg}
	\lbal
	u\in L^\infty((0,\infty);L^1(\Om)) \cap L^\frac{n+2}{n}_{loc}(\bom\times [0,\infty)) 
		\cap L^\frac{n+2}{n+1}_{loc}([0,\infty);W^{1,\frac{n+2}{n+1}}(\Om))
	\qquad \mbox{and} \\[1mm]
	v\in L^\infty(\Om\times (0,\infty)) \cap L^2_{loc}([0,\infty);W^{2,2}(\Om))
	\ear
  \ee
  such that $(u,v)$ forms a global weak solution of (\ref{0}) in the sense of Definition \ref{dw}.
  This solution has the additional property that there exists $T>0$ fulfilling  
  \be{15.1}
	u\in C^{2,1}(\bom\times [T,\infty))
	\qquad \mbox{and} \qquad
	v\in C^{2,1}(\bom\times [T,\infty)),  
  \ee
  and we have
  \be{15.2}
	u(\cdot,t)\to \mu:=\frac{1}{|\Om|} \io u_0
	\quad \mbox{in } L^\infty(\Om)
	\qquad \mbox{as } t\to\infty
  \ee
  and 
  \be{15.3}
	v(\cdot,t)\to 0
	\quad \mbox{in } L^\infty(\Om)
	\qquad \mbox{as } t\to\infty.
  \ee
\end{theo}
The convergence asserted by \eqref{15.2} and \eqref{15.3} in fact occurs at an exponential rate, see the discussion near \cite[Theorem 5]{lan_win_depleting} based on lower-dimensional precedents (\cite{zhang_li}, \cite{jiang}).
\mysection{Preliminaries}
\newcommand{\une}{u_{0\eps}}
\newcommand{\vne}{v_{0\eps}}
\newcommand{\defeq}{\colonequals}
Our analysis will operate on the regularized variants of (\ref{0}) given by
\be{0eps}
	\left\{ \begin{array}{ll}	
	u_{\eps t} =\Del\ueps - \na\cdot \Big( \frac{\ueps}{(1+\eps\ueps)^2} \na\veps\Big),
	\qquad & x\in\Om, \ t>0, \\[1mm]
	v_{\eps t} = \Del \veps - \frac{\ueps\veps}{1+\eps\ueps},
	\qquad & x\in\Om, \ t>0, \\[1mm]
	\frac{\pa \ueps}{\pa\nu}=\frac{\pa \veps}{\pa\nu}=0,
	\qquad & x\in\pO, \ t>0, \\[1mm]
	\ueps(x,0)=\une (x), \quad \veps(x,0)=\vne (x)
	\qquad & x\in\Om,
	\end{array} \right.
\ee 
for $\eps\in (0,1)$, where we let $\une$ and $\vne$ be chosen according to the following.
\begin{lem}\label{lem:approx}
  Let $u_0$ and $v_0$ satisfy \eqref{init}. Then for each $\eps\in(0,1)$ there are functions $\une$ and $\vne$ such that 
  \begin{subequations}\label{init-approx-u}
  \begin{align}
  	& \une\in C^\infty(\bom), \qquad \une\ge 0, \qquad \une \not\equiv 0 \qquad &&\text{for each } \eps\in(0,1),\label{init-u-regularity}\\
  	& \une \to u_0 \quad \text{ in } \Lom1 \quad &&\text{as } \eps\searrow 0,\label{init-u-L1}\\
  	& \io \une = \io u_0 \qquad &&\text{for each } \eps\in(0,1),\label{init-u-mass}\\
  	& \io \une \ln \une \le 1 + \io u_0\ln u_0  \qquad &&\text{for each } \eps\in(0,1),\label{init-u-LlogL}
  \end{align}
  \end{subequations}
  and that 
  \begin{subequations}\label{init-approx-v}
  \begin{align}
   	& \vne \in  C^\infty(\bom), \qquad \vne>0  \qquad &&\text{for each } \eps\in(0,1),\label{init-v-regularity}\\
   	& \vne \to v_0 \quad &&\text{in } \Lom1 \quad \text{as } \eps\searrow 0,\label{init-v-L1}\\
   	& \norm[\Lom p]{\vne}\le \norm[\Lom p]{v_0}&&\text{for each $p\in [1,\infty]$ and } \eps\in(0,1),\label{init-v-Lp-bounds}\\
   	& \io \f{|\na \vne|^2}{\vne}\le \io \f{|\na v_0|^2}{v_0} \qquad &&\text{for each } \eps\in(0,1).\label{init-v-nablaroot}
  \end{align}
  \end{subequations}
\end{lem}
\proof
 By $\eta_{\eps}$ we denote the usual mollifier and with $\alpha>0$ sufficiently small (so as to allow for the summand $1$ in \eqref{init-u-LlogL} in the end) we introduce $\tilde{u}_{0\eps}\defeq \eta_{\alpha\eps}\star u_0$. We then let $\une\defeq \f{\io u_0}{\io \tilde{u}_{0\eps}} \tilde{u}_{0\eps}$, ensuring \eqref{init-u-mass} along with regularity and approximation as required by \eqref{init-u-regularity}, \eqref{init-u-L1} and \eqref{init-u-LlogL} (cf. \cite[Thms.~2.29, 8.21]{adams}). 
  Regarding the approximation of $\ve$, we let $w_{\eps}\defeq \eta_{\alpha\eps}\star \sqrt{v_0}$ and define $\vne=w_{\eps}^2$, so that \eqref{init-approx-v} follows from the well-known properties of mollifications (eg., \cite[Thm.~2.29 (a) and (c)]{adams}). \qed
\qed
The approximate system \eqref{0eps} is chosen such that it has global classical solutions:

\begin{lem}\label{lem1}
  Let $\eps\in (0,1)$.
  For each $\eps\in (0,1)$, there exist
  \bas
	\lbal
	\ueps\in C^0(\bom\times [0,\infty)) \cap C^{2,1}(\bom\times (0,\infty))
	\qquad \mbox{and} \\[1mm]
	v\in \bigcap_{q>n} C^0([0,\infty);W^{1,q}(\Om)) \cap C^{2,1}(\bom\times (0,\infty))
	\ear
  \eas
  such that $\ueps\ge 0$ and $\veps>0$ in $\bom\times [0,\infty)$, and that (\ref{0eps}) is satisfied in the classical sense.
\end{lem}
\proof
  With the choices $S(u)=\f1{(1+\eps u)^2}$, $f\equiv 0$, $g(u,v)=v-\f{uv}{1+\eps u}$ taken there, the general chemotaxis existence 
  result of \cite[Lemma 3.1]{BBTW} yields local existence of solutions $(\ue,\ve)$ of the desired regularity 
  on $\Om\times(0,\Tmaxe)$ where, for $q>n$, either $\Tmaxe=\infty$ or 
  $\lim_{t\upto\Tmaxe} \kl{\norm[\Lom\infty]{\ue(\cdot,t)} +\norm[\Wom1q]{\ve(\cdot,t)}}=\infty$. 
  For each fixed $\eps\in(0,1)$, $-\f{\ue\ve}{1+\eps\ue}$ is a bounded function on $\Ombar\times[0,\Tmaxe)$, 
  hence $\norm[\Wom1q]{\ve(\cdot,t)}$ remains bounded. 
  From the resulting boundedness of $\f{\ue}{(1+\eps\ue)^2}\na \ve$, semigroup estimates similarly give rise to boundedness of 
  $t\mapsto \norm[\Lom\infty]{\ue(\cdot,t)}$, so that $\Tmaxe=\infty$ follows.
\qed
Throughout the sequel, without explicit mentioning we let $(\ueps)_{\eps\in (0,1)}$ and $(\veps)_{\eps\in (0,1)}$ be as
provided by Lemma \ref{lem1}. Some first boundedness information is obtained quickly.
\begin{lem}\label{lem2}
  Let $\eps\in (0,1)$ and $p\in[1,\infty]$. Then
  \be{mass}
	\io \ueps(\cdot,t) = \io u_0 \equiv \mu |\Om|
	\qquad \mbox{for all } t>0
  \ee
  and
  \be{vinfty}
	\|\veps(\cdot,t)\|_{L^p(\Om)} \le \|\veps(\cdot,t_0)\|_{L^p(\Om)} \le \|v_0\|_{L^p(\Om)}
	\qquad \mbox{for all $t_0\ge 0$ and $t>t_0$,}
  \ee
  and we have
  \be{uv}
	\int_0^\infty \io \frac{\ueps\veps}{1+\eps\ueps} \le \io v_0
  \ee
  as well as
  \be{gradv}
	\int_0^\infty \io |\na \veps|^2 \le \frac{1}{2} \io v_0^2.
  \ee
\end{lem}
\proof
  Integrating the first equation of \eqref{0eps} together with \eqref{init-u-mass} yields \eqref{mass}, while testing the second 
  equation by $1$ or $\ve$ or $\ve^{p-1}$ by \eqref{init-v-Lp-bounds} shows \eqref{uv}, \eqref{gradv} and \eqref{vinfty} 
  for $p<\infty$, respectively. For $p=\infty$, \eqref{vinfty} results from a comparison principle, again together with 
  \eqref{init-v-Lp-bounds}.
\qed
\mysection{Refined energy analysis}\label{sec:energy}
Our energy analysis will be launched by the following basic observation.
\begin{lem}\label{lem3}
  For any choice of $\eps\in (0,1)$, we have
  \be{en}
	\Feps'(t) + \D_\eps(t) = \frac{1}{2} \int_{\pO} \frac{1}{\veps} \frac{\pa |\na\veps|^2}{\pa\nu}
	\qquad \mbox{for all } t>0,
  \ee
  where
  \be{F}
	\Feps(t):=\io \ueps(\cdot,t) \ln \frac{\ueps(\cdot,t)}{\mu} 
	+ \frac{1}{2} \io \frac{|\na \veps(\cdot,t)|^2}{\veps(\cdot,t)},
	\qquad t\ge 0,
  \ee
  and
  \be{D}
	\D_\eps(t):=\io \frac{|\na\ueps(\cdot,t)|^2}{\ueps(\cdot,t)}
	+ \io \veps(\cdot,t) |D^2 \ln\veps(\cdot,t)|^2
	+ \frac{1}{2} \io \frac{\ueps(\cdot,t)}{1+\eps\ueps(\cdot,t)} \frac{|\na \veps(\cdot,t)|^2}{\veps(\cdot,t)},
	\qquad t>0.
  \ee
\end{lem}
\proof
  This can be seen by straightforward computation on the basis of (\ref{0eps}) (cf.~\cite[Lemma 3.2]{win_CPDE} for details).
\qed
The dissipative term involving $D^2\ln \ve$ in \eqref{D} already controls more directly useful expressions:
\begin{lem}\label{lem33}
  Let $\vp\in C^2(\bom)$ be such that $\vp>0$ in $\bom$ and $\frac{\pa\vp}{\pa\nu}=0$ on $\pO$. Then
  \be{33.1}
	\io \frac{|\na\vp|^4}{\vp^3} \le (2+\sqrt{n})^2 \io \vp |D^2\ln\vp|^2
  \ee
  and
  \be{33.2}
	\io \frac{|D^2\vp|^2}{\vp} \le (2n+8\sqrt{n}+10) \io \vp|D^2\ln\vp|^2.
  \ee
\end{lem}
\proof
  The inequality in (\ref{33.1}) has been derived in \cite[Lemma 3.3]{win_CPDE}.
  For a verification of (\ref{33.2}) on the basis of this, we only need to apply the inequality $|a-b|^2 \ge \frac{1}{2}a^2-b^2$,
  valid for all $(a,b)\in\R^2$, to see that
  \bas
	\io \vp |D^2\ln\vp|^2
	&=& \sum_{i,j=1}^n \io \vp \Big| \frac{\pa_{x_i x_j} \vp}{\vp} - \frac{\pa_{x_i} \vp \pa_{x_j} \vp}{\vp^2} \Big|^2 \\
	&\ge& \frac{1}{2} \sum_{i,j=1}^n \io \frac{|\pa_{x_i x_j} \vp|^2}{\vp}
	- \sum_{i,j=1}^n \io \frac{|\pa_{x_i}\vp|^2 |\pa_{x_j}\vp|^2}{\vp^3} \\
	&=& \frac{1}{2} \io \frac{|D^2\vp|^2}{\vp} - \io \frac{|\na\vp|^4}{\vp^3}.
  \eas
  Therefore, (\ref{33.2}) indeed results from (\ref{33.1}), because $2+ 2\cdot(2+\sqrt{n})^2=2n+8\sqrt{n}+10$.
\qed
Now the possibly most crucial observation for our large time analysis is contained in the following simple but useful 
statement on one-sided control that can be achieved for boundary integrals of the form in (\ref{en}).
\begin{lem}\label{lem44}
  Let $\eta>0$. Then there exists $C(\eta)>0$ such that whenever $\vp\in C^2(\bom)$ is positive in $\bom$ with
  $\frac{\pa\vp}{\pa\nu}=0$ on $\pO$,
  \be{44.1}
	\int_{\pO} \frac{1}{\vp} \frac{\pa |\na\vp|^2}{\pa\nu}
	\le \eta \io \frac{|D^2\vp|^2}{\vp}
	+ \eta \io \frac{|\na\vp|^4}{\vp^3}
	+ C(\eta) \io |\na\vp|.	
  \ee
\end{lem}
\proof
  According to a boundary trace embedding inequality, we can fix $c_1>0$ such that
  \be{44.3}
	\int_{\pO} |\psi| \le c_1 \io |\na\psi| + c_1 \io |\psi|
	\qquad \mbox{for all } \psi\in C^1(\bom),
  \ee
  and from \cite{mizoguchi_souplet} we obtain $c_2>0$ fulfilling
  \be{44.4}
	\frac{\pa |\na\vp|^2}{\pa\nu} \le c_2 |\na\vp|^2
	\quad \mbox{on } \pO
	\qquad \mbox{for all $\vp\in C^2(\bom)$ such that $\frac{\pa\vp}{\pa\nu}=0$ on } \pO.
  \ee
  Fixing any such $\vp$ which is positive in $\Ombar$, by combination of (\ref{44.4}) with (\ref{44.3}) we see that
  \bea{44.5}
	\int_{\pO} \frac{1}{\vp} \frac{\pa |\na\vp|^2}{\pa\nu}
	&\le& c_2 \int_{\pO} \frac{|\na\vp|^2}{\vp} \nn\\
	&\le& c_1 c_2 \io \Big| \na \frac{|\na\vp|^2}{\vp}\Big|
	+ c_1 c_2 \io \frac{|\na\vp|^2}{\vp} \nn\\
	&=& c_1 c_2 \io \Big| \frac{2}{\vp} D^2\vp\cdot\na\vp - \frac{1}{\vp^2} |\na\vp|^2 \na\vp\Big|
	+ c_1 c_2 \io \frac{|\na\vp|^2}{\vp} \nn\\
	&\le& 2c_1 c_2 \io \frac{|D^2\vp| \cdot |\na\vp|}{\vp}
	+ c_1 c_2 \io \frac{|\na\vp|^3}{\vp^2}
	+ c_1 c_2 \io \frac{|\na\vp|^2}{\vp}.
  \eea
  Given $\eta>0$, by means of Young's inequality we see that here
  \bas
	2c_1 c_2 \io \frac{|D^2\vp| \cdot |\na\vp|}{\vp}
	\le \eta\io \frac{|D^2\vp|^2}{\vp}
	+ \frac{c_1^2 c_2^2}{\eta} \io \frac{|\na\vp|^2}{\vp}
  \eas
  and
  \bas
	c_1 c_2 \io \frac{|\na\vp|^3}{\vp^2}
	\le \frac{\eta}{2} \io \frac{|\na\vp|^4}{\vp^3}
	+ \frac{c_1^2 c_2^2}{2\eta} \io \frac{|\na\vp|^2}{\vp},
  \eas
  and that writing $c_3(\eta):=\frac{3 c_1^2 c_2^2}{2\eta} + c_1c_2$ we have
  \bas
	c_3(\eta) \io \frac{|\na\vp|^2}{\vp}
	&=& \io \Big\{ \frac{\eta}{2} \frac{|\na\vp|^4}{\vp^3} \Big\}^\frac{1}{3} \cdot 
		c_3(\eta) \cdot \Big(\frac{2}{\eta}\Big)^\frac{1}{3} \cdot |\na\vp|^\frac{2}{3} \\
	&\le& \frac{\eta}{2} \io \frac{|\na\vp|^4}{\vp^3}
	+ \Big\{ c_3(\eta) \cdot \Big(\frac{2}{\eta}\Big)^\frac{1}{3} \Big\}^\frac{3}{2} \io |\na\vp|.
  \eas
  Therefore, (\ref{44.5}) implies (\ref{44.1}) if we let
  $C(\eta):=\sqrt{\frac{2c_3^3(\eta)}{\eta}}$.
\qed
When combined with a simple argument that conveniently estimates the first dissipated integral in (\ref{en})-(\ref{D})
from below, Lemma \ref{lem44} entails the following.
\begin{lem}\label{lem4}
  There exist $\gamma>0$ and $\Gamma>0$ such that with $(\Feps)_{\eps\in (0,1)}$ taken from (\ref{F}) and $\mu$ from \eqref{mass},
  \bea{4.1}
	& & \hs{-20mm}
	\Feps'(t)
	+ \gamma \io \frac{|\na\ueps|^2}{\ueps}
	+ \gamma \cdot\bigg\{ \io |\ueps-\mu| \bigg\}^2
	+ \gamma \io \frac{|D^2\veps|^2}{\veps}
	+ \gamma \io \frac{\ueps}{1+\eps\ueps} \frac{|\na\veps|^2}{\veps} \nn\\
	&\le& \Gamma \io |\na\veps|
	\qquad \mbox{for all $t>0$ and } \eps\in (0,1).
  \eea
\end{lem}
\proof
  The combination of a Poincar\'e inequality with (\ref{mass}) yields $c_1>0$ such that
  \bas
	\io |\ueps-\mu| \le c_1 \io |\na\ueps|
	\qquad \mbox{for all $t>0$ and } \eps\in (0,1),
  \eas
  whence by the Cauchy-Schwarz inequality, and again by (\ref{mass}),
  \bas
	\bigg\{ \io |\ueps-\mu| \bigg\}^2
	&\le& c_1^2 \cdot \bigg\{ \io |\na\ueps| \bigg\}^2 \\
	&\le& c_1^2 \cdot \bigg\{ \io \ueps \bigg\} \cdot \io \frac{|\na\ueps|^2}{\ueps} \\
	&=& c_1^2 \mu |\Om| \io \frac{|\na\ueps|^2}{\ueps}
	\qquad \mbox{for all $t>0$ and } \eps\in (0,1).
  \eas
  From Lemma \ref{lem3} and Lemma \ref{lem33} we therefore readily obtain that with some $c_2>0$,
  \bas
	& & \hs{-20mm}
	\Feps'(t)
	+ \frac{1}{2} \io \frac{|\na\ueps|^2}{\ueps}
	+ c_2 \cdot \bigg\{ \io |\ueps-\mu| \bigg\}^2
	+ c_2 \io \frac{|D^2\veps|^2}{\veps}
	+ c_2 \io \frac{|\na\veps|^4}{\veps^3}
	+ c_2 \io \frac{\ueps}{1+\eps\ueps} \frac{|\na\veps|^2}{\veps} \\
	&\le& \frac{1}{2} \int_{\pO} \frac{1}{\veps} \frac{\pa |\na\veps|^2}{\pa\nu}
	\qquad \mbox{for all $t>0$ and } \eps\in (0,1),
  \eas
  whence an application of Lemma \ref{lem44} to $\eta:=c_2$ already leads to (\ref{4.1}) with
  $\gamma:=\min\{ \frac{1}{2} \, , \, \frac{c_2}{2}\}$ and some suitably large $\Gamma>0$.
\qed
As \eqref{gradv} asserts some large time decay property of the right-hand side in \eqref{4.1},
an integration in time provides bounds which depend on the length of the integration interval in a sublinear manner.
\begin{lem}\label{lem5}
  There exist $C>0$ such that
  \be{5.1}
	\int_0^T \io |\ueps-\mu| \le CT^\frac{3}{4}
	\qquad \mbox{for all $T>1$ and } \eps\in (0,1)
  \ee
  and
  \be{5.2}
	\int_0^T \io \frac{|\na\ueps|^2}{\ueps} \le C\sqrt{T}
	\qquad \mbox{for all $T>1$ and } \eps\in (0,1)
  \ee
  as well as
  \be{5.4}
	\int_0^T \io |D^2\veps|^2 \le C\sqrt{T}
	\qquad \mbox{for all $T>1$ and } \eps\in (0,1)
  \ee
  and
  \be{5.3}
	\int_0^T \io \frac{\ueps}{1+\eps\ueps} \frac{|\na\veps|^2}{\veps} \le C\sqrt{T}
	\qquad \mbox{for all $T>1$ and } \eps\in (0,1).
  \ee
\end{lem}
\proof
  With $\gamma>0$ and $\Gamma>0$ as given by Lemma \ref{lem4}, for $\eps\in (0,1)$ and $t>0$ we let
  \bas
	\geps(t):=\gamma \io \frac{|\na\ueps(\cdot,t)|^2}{\ueps(\cdot,t)}
	+ \gamma\cdot\bigg\{ \io |\ueps(\cdot,t)-\mu| \bigg\}^2
	+ \gamma \io \frac{|D^2\veps(\cdot,t)|^2}{\veps(\cdot,t)} 
	+ \gamma \io \frac{\ueps(\cdot,t)}{1+\eps\ueps(\cdot,t)} \frac{|\na\veps(\cdot,t)|^2}{\veps(\cdot,t)}
  \eas
  and
  \bas
	\heps(t):=\Gamma \io |\na \veps(\cdot,t)|,
  \eas
  and note that according to the Cauchy-Schwarz inequality and (\ref{gradv}),
  \bas
	\int_0^T \heps(t) dt 
	&\le& \Gamma \sqrt{|\Om|} \sqrt{T} \cdot \bigg\{ \int_0^T \io |\na\veps|^2 \bigg\}^\frac{1}{2} \\[2mm]
	&\le& c_1\sqrt{T}
	\qquad \mbox{for all $T>0$ and } \eps\in (0,1)
  \eas
  with $c_1:=\Gamma\sqrt{|\Om|} \cdot \Big\{ \frac{1}{2} \io v_0^2\Big\}^\frac{1}{2}$.
  An integration of (\ref{4.1}) thus shows that by \eqref{init-u-LlogL} and \eqref{init-v-nablaroot}
  \bas
	\Feps(T) + \int_0^T \geps(t) dt
	&\le& \Feps(0) + \int_0^T \heps(t) dt \\
	&\le& 1+\io u_0\ln\frac{u_0}{\mu} + \frac{1}{2} \io \frac{|\na v_0|^2}{v_0}
	+ c_1\sqrt{T}
	\qquad \mbox{for all $T>0$ and } \eps\in (0,1),
  \eas
  so that all claimed inequalities result upon observing that $\Feps$ is nonnegative by (\ref{mass}), estimating
  $\frac{1}{\veps} \ge \frac{1}{\|v_0\|_{L^\infty(\Om)}}$ by (\ref{vinfty}), and using that
  \bas
	\int_0^T \io |\ueps-\mu|
	\le \sqrt{T} \int_0^T \bigg\{ \io |\ueps-\mu|\bigg\}^2
	\le \sqrt{\frac{T}{\gamma}} \cdot\bigg\{ \int_0^T \geps(t) dt \bigg\}^\frac{1}{2}
	\qquad \mbox{for all $T>0$ and } \eps\in (0,1)
  \eas
  thanks to the Cauchy-Schwarz inequality.
\qed
\mysection{Passing to the limit. Global existence}
In this section aiming at the mere construction of global solutions, we may ignore the particular dependence of the estimates in
Section \ref{sec:energy} on $T$, and collect the above findings to infer the following set of bounds. Notable among these is \eqref{6.3}, which -- in contrast to the treatment of this term in \cite{win_CPDE} -- relies on the dissipative term controlled in \eqref{5.3} and not on boundedness of the $\ue$ and the $\nabla\ve$ obtained separately. 
This allows to achieve some $L^1$ compactness property regardless of the spatial dimension.
\begin{lem}\label{lem6}
  For all $T>1$ there exists $C(T)>0$ with the property that
  \be{6.1}
	\int_0^T \io \ueps^\frac{n+2}{n} \le C(T)
	\qquad \mbox{for all } \eps\in (0,1),
  \ee
  that
  \be{6.2}
	\int_0^T \io |\na\ueps|^\frac{n+2}{n+1} \le C(T)
	\qquad \mbox{for all } \eps\in (0,1),
  \ee
  and that
  \be{6.3}
	\int_0^T\io \Big| \frac{\ueps}{(1+\eps\ueps)^2} \na\veps\Big|^\frac{n+2}{n+1} \le C(T)
	\qquad \mbox{for all } \eps\in (0,1).
  \ee
\end{lem}
\proof
  Since from the Gagliardo-Nirenberg inequality we obtain $c_1>0$ such that
  \bas
	\io \ueps^\frac{n+2}{n}
	= \|\sqrt{\ueps}\|_{L^\frac{2(n+2)}{n}(\Om)}^\frac{2(n+2)}{n}
	\le c_1 \|\na\sqrt{\ueps}\|_{L^2(\Om)}^2 \|\sqrt{\ueps}\|_{L^2(\Om)}^\frac{4}{n}
	+ c_1\|\sqrt{\ueps}\|_{L^2(\Om)}^\frac{2(n+2)}{n}
	\quad \mbox{for all $t>0$ and } \eps\in (0,1),
  \eas
  the inequality in (\ref{6.1}) follows from (\ref{5.2}) and (\ref{mass}) after integrating here.
  We thereupon invoke Young's inequality to see that
  \bas
	\int_0^T \io \Big|\frac{\ueps}{(1+\eps\ueps)^2} \na\veps\Big|^\frac{n+2}{n+1}
	&=& \int_0^T \io \Big\{ \frac{\ueps}{1+\eps\ueps} |\na\veps|^2 \Big\}^\frac{n+2}{2(n+1)} 
		\cdot \Big\{ \frac{\ueps}{(1+\eps\ueps)^3} \Big\}^\frac{n+2}{2(n+1)} \\
	&\le& \int_0^T \io \frac{\ueps}{1+\eps\ueps} |\na\veps|^2
	+ \int_0^T \io \Big\{ \frac{\ueps}{(1+\eps\ueps)^3} \Big\}^\frac{n+2}{n} \\
	&\le& \int_0^T \io \frac{\ueps}{1+\eps\ueps} |\na\veps|^2
	+ \int_0^T \io \ueps^\frac{n+2}{n}
	\qquad \mbox{for all $T>1$ and } \eps\in (0,1),
  \eas
  so that (\ref{6.3}) results from (\ref{5.3}) and (\ref{6.1}).
  The estimate in (\ref{6.2}) can be derived in much the same manner on the basis of (\ref{5.2}) and (\ref{6.1}).
\qed
In preparation of an Aubin--Lions type compactness argument, we next provide some bounds for temporal derivatives.
\begin{lem}\label{lem7}
  Let $T>1$. Then there exists $C(T)>0$ such that
  \be{7.1}
	\int_0^T \big\|u_{\eps t}(\cdot,t)\big\|_{(W^{1,n+2}(\Om))^\star}^\frac{n+2}{n+1} dt \le C(T)
	\qquad \mbox{for all } \eps\in (0,1)
  \ee
  and
  \be{7.2}
	\int_0^T \io |v_{\eps t}|^\frac{n+2}{n+1} \le C(T)
	\qquad \mbox{for all } \eps\in (0,1).
  \ee
\end{lem}
\proof
  For the derivation of (\ref{7.1}), we only need to observe that in the identity $u_{\eps t}=\na\cdot a_\eps(x,t)$, $\eps\in (0,1)$,
  according to Lemma \ref{lem6} the functions $a_\eps:=\na\ueps-\frac{\ueps}{(1+\eps\ueps)^2} \na\veps$, $\eps\in (0,1)$,
  have the property that for each $T>0$ there exists $c_1(T)>0$ such that
  \bas
	\int_0^T \io |a_\eps|^\frac{n+2}{n+1} \le c_1(T)
	\qquad \mbox{for all } \eps\in (0,1).
  \eas
  Therefore, (\ref{7.1}) follows in a standard manner from the definition of norms in dual spaces, 
  while (\ref{7.2}) directly results from the second equation in (\ref{0eps}) and the fact that both
  $(\Del\veps)_{\eps\in (0,1)}$ and $\big(\frac{\ueps\veps}{1+\eps\ueps})_{\eps\in (0,1)}$
  are bounded in $L^\frac{n+2}{n+1}(\Om\times (0,T))$ for all $T>1$ thanks to (\ref{5.4}), (\ref{6.1}), (\ref{vinfty})
  and the inequality $\frac{n+2}{n+1} \le 2$.
\qed
With these prerequisites at hand, we can derive the following by a suitable limit passage in (\ref{0eps}) in a straightforward
manner.
\begin{lem}\label{lem9}
  There exist $(\eps_j)_{j\in\N}\subset (0,1)$ as well as nonnegative functions $u$ and $v$ which satisfy (\ref{reg}) 
  as well as
  \bas
	u\na v\in L^\frac{n+2}{n+1}_{loc}(\bom\times [0,\infty);\R^n),
  \eas
  such that $\eps_j\searrow 0$ as $j\to\infty$, and that
  \begin{eqnarray}
	& & \ueps\to u
	\qquad \mbox{in $L^1_{loc}(\bom\times [0,\infty))$ and a.e.~in } \Om\times (0,\infty),
	\label{9.1} \\
	& & \na\ueps\wto\na u
	\qquad \mbox{in $L^1_{loc}(\bom\times [0,\infty))$},
	\label{9.2} \\
	& & \veps\to v
	\qquad \mbox{a.e.~in } \Om\times (0,\infty),
	\label{9.3} \\
	& & \veps \wsto v
	\qquad \mbox{in $L^\infty(\Om\times (0,\infty))$,}
	\label{9.4} \\
	& & \na\veps\to\na v
	\qquad \mbox{in $L^1_{loc}(\bom\times [0,\infty))$ and a.e.~in } \Om\times (0,\infty),
	\label{9.5} \\
	& & \frac{\ueps}{(1+\eps\ueps)^2} \na\veps \to u\na v
	\qquad \mbox{in $L^1_{loc}(\bom\times [0,\infty))$}
	\label{9.6}
  \end{eqnarray}
  as $\eps=\eps_j\searrow 0$. 
  Moreover, $(u,v)$ is a global weak solution of (\ref{0}) in the sense of Definition \ref{dw}.
\end{lem}
\proof
  Since (\ref{6.2}), (\ref{mass}) and (\ref{7.1}) assert boundedness of $(\ueps)_{\eps\in (0,1)}$ in
  $L^\frac{n+2}{n+1}((0,T);W^{1,\frac{n+2}{n+1}}(\Om))$ and of $(u_{\eps t})_{\eps\in (0,1)}$ in
  $L^\frac{n+2}{n+1}((0,T);(W^{1,n+2}(\Om))^\star)$ for all $T>1$,
  and since (\ref{5.4}), (\ref{vinfty}) and (\ref{7.2}) ensure boundedness of $(\veps)_{\eps\in (0,1)}$
  in $L^2((0,T);W^{2,2}(\Om))$ and of $(v_{\eps t})_{\eps\in (0,1)}$ in $L^\frac{n+2}{n+1}(\Om\times (0,T))$ for each $T>1$,
  two applications of an Aubin-Lions lemma (\cite{simon})
  yield $(\eps_j)_{j\in\N}\subset (0,1)$, $u\in L^\frac{n+2}{n+1}_{loc}([0,\infty);W^{1,\frac{n+2}{n+1}}(\Om))$
  and $v\in L^2_{loc}([0,\infty);W^{2,2}(\Om))$ such that $\eps_j\searrow 0$ as $j\to\infty$, and that
  (\ref{9.1}), (\ref{9.2}), (\ref{9.3}), and (\ref{9.5}) hold as $\eps=\eps_j\searrow 0$, whence
  clearly $u$ and $v$ are nonnegative a.e.~in $\Om\times (0,\infty)$.\abs
  Since we thus also know that, in particular, $\frac{\ueps}{(1+\eps\ueps)^2} \na\veps \to u\na v$ a.e.~in $\Om\times (0,\infty)$
  as $\eps=\eps_j\searrow 0$, and since $\big(\frac{\ueps}{(1+\eps\ueps)^2} \na\veps\big)_{\eps\in (0,1)}$
  is uniformly integrable over $\Om\times (0,T)$ for all $T>1$ according to (\ref{6.3}), we obtain (\ref{9.6}) as a consequence
  of the Vitali convrgence theorem, whereas (\ref{9.4}) similarly results from (\ref{9.3}) and (\ref{vinfty}).\abs
  Deriving (\ref{wu}) and (\ref{wv}) from the corresponding weak formulations of the respective sub-problems of (\ref{0eps}) 
  can thereupon be accomplished on the basis of \eqref{init-approx-u}, \eqref{init-approx-v}, (\ref{9.1}), (\ref{9.2}) and (\ref{9.4})-(\ref{9.6}) in a straightforward fashion.
\qed
\mysection{Large time behavior}
As a first step regarding the long term behavior, we observe that the $\Lom1$ norm of $\ve$ becomes small, uniformly with respect to $\eps$. Here sublinearity of the growth observed in Lemma~\ref{lem5} is important.
\begin{lem}\label{lem11}
  Let $\del>0$. Then there exist $T_\star(\del)>0$ and $\eps_\star(\del)\in (0,1)$ such that
  \be{11.1}
	\io \veps\big(\cdot,T_\star(\del)\big) \le \del
	\qquad \mbox{for all } \eps\in (\eps_j)_{j\in\N} \cap (0,\eps_\star(\del)).
  \ee
\end{lem}
\proof
  We abbreviate $c_1:=\io v_0$ and $c_2:=\|v_0\|_{L^\infty(\Om)}$, and using Lemma \ref{lem5} we find $c_3>0$ such that
  \be{11.2}
	\int_0^T \io |\ueps-\mu| \le c_3 T^\frac{3}{4}
	\qquad \mbox{for all $T>1$ and } \eps\in (0,1).
  \ee
  Given $\delta>0$, we then rely on our overall assumption that $u_0$ be nontrivial in choosing 
  $T_\star=T_\star(\del)>1$ large enough such that
  \be{11.3}
	\frac{c_1}{T_\star} + \frac{c_2 c_3}{T_\star^\frac{1}{4}} \le \frac{\mu\del}{2},
  \ee
  and in fixing some small $\eta=\eta(\del)>0$ in such a way that also
  \be{11.4}
	\eta\le\frac{\mu\del}{2}.
  \ee
  Since
  \bas
	\frac{\eps\ueps^2 \veps}{1+\eps\ueps} \le \ueps\veps \le c_2 \ueps
	\quad \mbox{in } \Om\times (0,\infty)
	\qquad \mbox{for all } \eps\in (0,1)
  \eas
  due to (\ref{vinfty}), using that $\ueps\to u$ in $L^1_{loc}(\bom\times [0,\infty))$ as $\eps=\eps_j\searrow 0$ by 
  Lemma \ref{lem9} we may employ a version of the dominated convergence theorem to see that
  \bas
	\int_0^{T_\star} \io \frac{\eps\ueps^2 \veps}{1+\eps\ueps} \to 0
	\qquad \mbox{as } \eps=\eps_j\searrow 0,
  \eas
  whence we can finally pick $\eps_\star=\eps_\star(\del)\in (0,1)$ suitably small fulfilling
  \be{11.5}
	\frac{1}{T_\star} \int_0^{T_\star} \io \Big( \ueps-\frac{\ueps}{1-\eps\ueps}\Big) \veps
	= \frac{1}{T_\star} \int_0^{T_\star} \io \frac{\eps\ueps^2 \veps}{1+\eps\ueps}
	\le \eta
	\qquad \mbox{for all } \eps\in (\eps_j)_{j\in\N} \cap (0,\eps_\star).
  \ee
  Now decomposing
  \bas
	\frac{1}{T_\star} \int_0^{T_\star} \io \frac{\ueps\veps}{1+\eps\ueps}
	&=& \frac{\mu}{T_\star} \int_0^{T_\star} \io \veps \\
	& & + \frac{1}{T_\star} \int_0^{T_\star} \io (\ueps-\mu) \veps 
	- \frac{1}{T_\star} \int_0^{T_\star} \io \Big(\ueps-\frac{\ueps}{1+\eps\ueps}\Big) \veps
	\qquad \mbox{for } \eps\in (0,1),
  \eas
  from (\ref{uv}), (\ref{11.2}) and (\ref{11.5}) we obtain that
  \bas
	\frac{\mu}{T_\star} \int_0^{T_\star} \io \veps
	&\le& \frac{1}{T_\star} \io v_0
	- \frac{1}{T_\star} \int_0^{T_\star} \io (\ueps-\mu) \veps
	+ \frac{1}{T_\star} \int_0^{T_\star} \io \Big(\ueps-\frac{\ueps}{1+\eps\ueps}\Big)\veps \\
	&\le& \frac{c_1}{T_\star}
	+ \frac{c_2}{T_\star} \int_0^{T_\star} \io |\ueps-\mu|
	+ \frac{1}{T_\star} \int_0^{T_\star} \io \Big(\ueps-\frac{\ueps}{1+\eps\ueps}\Big)\veps \\
	&\le& \frac{c_1}{T_\star}
	+ \frac{c_2 c_3}{T_\star^\frac{1}{4}}
	+ \eta
	\qquad \mbox{for all } \eps\in (\eps_j)_{j\in\N} \cap (0,\eps_\star),
  \eas
  whence (\ref{11.3}) an (\ref{11.4}) imply that
  \bas
	\frac{1}{T_\star} \int_0^{T_\star} \io \veps
	\le \frac{1}{\mu} \cdot \Big\{ \frac{\mu\del}{2} + \frac{\mu\del}{2}\Big\} =\del
	\qquad \mbox{for all } \eps\in (\eps_j)_{j\in\N} \cap (0,\eps_\star).
  \eas
  For any $\eps\in (\eps_j)_{j\in\N}\cap (0,\eps_\star)$, we can therefore find $t_0(\eps)\in (0,T_\star)$ such that
  \bas
	\io \veps \big(\cdot,t_0(\eps)\big) \le \del,
  \eas
  so that (\ref{11.1}) follows upon recalling that for any such $\eps$ we have
  $\io \veps(\cdot,t) \le \io \veps(\cdot,t_0(\eps))$ for all $t>t_0(\eps)$ by (\ref{vinfty}).
\qed
Parabolic smoothing carries the above statement over to $L^\infty$ topologies:
\begin{lem}\label{lem12}
  For each $\del>0$ there exist $T_{\star\star}(\del)>0$ and $\eps_{\star\star}(\del)\in (0,1)$ such that
  \be{12.1}
	\|\veps(\cdot,t)\|_{L^\infty(\Om)} \le \del
	\qquad \mbox{for all $t>T_{\star\star}(\del)$ and } \eps\in (\eps_j)_{j\in\N} \cap (0,\eps_{\star\star}(\del)).
  \ee
\end{lem}
\proof
  By means of a known smoothing property of the Neumann heat semigroup $(e^{t\Del})_{t\ge 0}$ on $\Om$,
  we can pick $c_1>0$ such that
  \be{12.2}
	\|e^\Del \vp\|_{L^\infty(\Om)} \le c_1 \|\vp\|_{L^1(\Om)}
	\qquad \mbox{for all } \vp\in C^0(\bom).
  \ee
  If, for fixed $\del>0$, we thereupon let
  \bas
	T_{\star\star}\equiv T_{\star\star}(\del):=T_\star\Big(\frac{\del}{c_1}\Big) +1
	\qquad \mbox{and} \qquad
	\eps_{\star\star}\equiv \eps_{\star\star}(\del):=\eps_\star \Big(\frac{\del}{c_1}\Big),
  \eas
  with $\big(T_\star(\wh{\del})\big)_{\wh{\del}>0}$ and 
  $\big(\eps_\star(\wh{\del})\big)_{\wh{\del}>0}$ as given by Lemma \ref{lem11},
  then observing that $v_{\eps t} \le \Del\veps$ in $\Om\times (0,\infty)$ for all $\eps\in (0,1)$ we can use (\ref{11.2})
  to estimate
  \bas
	\veps(\cdot,T_{\star\star})
	&\le& e^\Del \veps(\cdot,T_{\star\star}-1) \\
	&\le& \|e^\Del \veps(\cdot,T_{\star\star}-1)\|_{L^\infty(\Om)} \\
	&\le& c_1 \|\veps(\cdot,T_{\star\star}-1)\|_{L^1(\Om)}
	\quad \mbox{in } \Om
	\qquad \mbox{for all } \eps\in (0,1).
  \eas
  As our choices of $T_{\star\star}$ and $\eps_{\star\star}$ ensure that in line with (\ref{11.1}) we have
  \bas
	\|\veps(\cdot,T_{\star\star}-1)\|_{L^1(\Om)}
	= \io \veps\Big( \cdot,T_\star\Big(\frac{\del}{c_1}\Big) \Big)
	\le \frac{\del}{c_1}
	\qquad \mbox{for all } \eps\in (\eps_j)_{j\in\N} \cap (0,\eps_{\star\star}),
  \eas
  this implies that
  \bas
	\veps(\cdot,T_{\star\star}) \le \del
	\quad \mbox{in } \Om
	\qquad \mbox{for all } \eps\in (\eps_j)_{j\in\N} \cap (0,\eps_{\star\star}),
  \eas
  and that thus
  \bas
	\|\veps(\cdot,T_{\star\star})\|_{L^\infty(\Om)} \le \del
	\qquad \mbox{for all } \eps\in (\eps_j)_{j\in\N} \cap (0,\eps_{\star\star}),
  \eas
  because $\veps\ge 0$ for all $\eps\in (0,1)$.
  In view of (\ref{vinfty}), this establishes (\ref{12.1}).
\qed
The upcoming boundedness result regarding $\Lom p$ norms of $\ue$ will rely on the following observation resulting from an ODE comparison. 
\begin{lem}\label{lem:ode}
  For every $a,b,\tau>0$ and $\lambda>1$ there is $C(a,b,\tau,\lambda)>0$ such that for every $T\in(0,\infty]$ every solution 
  $y\in C^0([0,T))\cap C^1((0,T))$   of 
  \begin{equation}\label{ode}
	y'+ay^{\lambda}\le b
  \end{equation}
  satisfies $y(t)\le C$ for all $t\in (\tau,T)$.
\end{lem}
\proof
  A proof based on explicit computations via separation of variables can be found detailed in \cite[Lemma 3.8]{knosalla_lankeit}.
\qed
\begin{lem}\label{lem13}
  For all $p>3$
  there exist $T_{\star\star\star}(p)>0$, $\eps_{\star\star\star}(p)\in (0,1)$ and $C(p)>0$ such that
  \be{13.1}
	\io \ueps^p(\cdot,t) \le C(p)
	\qquad \mbox{for all $t>T_{\star\star\star}(p)$ and } \eps\in (\eps_j)_{j\in\N} \cap (0,\eps_{\star\star\star}(p)),
  \ee
  and that moreover
  \be{13.2}
	\int_{T_{\star\star\star}(p)}^\infty \io (\ueps+1)^{p-2} |\na\ueps|^2 \le C(p)
	\qquad \mbox{for all } \eps\in (\eps_j)_{j\in\N} \cap (0,\eps_{\star\star\star}(p)).
  \ee
\end{lem}
\proof
  Given $p>3$, we use that $\f{(2p)^2}{3p(p-1)}<2$ in fixing $\del=\del(p)>0$ such that
  \be{13.4}
	\f{(2p+\del p(p-1))^2}{3p(p-1)}+p\del-2\le 0,
  \ee
  and invoke Lemma \ref{lem12} to choose some $T^{(1)}=T^{(1)}(p)>0$ and some $\eps^{(1)}=\eps^{(1)}(p)\in (0,1)$ fulfilling
  \be{13.5}
	\veps \le \frac{\del}{2}
	\quad \mbox{in } \Om\times (T^{(1)},\infty)
	\qquad \mbox{for all } \eps\in (\eps_j)_{j\in\N} \cap (0,\eps^{(1)}).
  \ee
  Then integrating by parts using (\ref{0eps}) we find that for all $t>T^{(1)}$ and $\eps\in (\eps_j)_{j\in\N} \cap (0,\eps^{(1)})$,
  \begin{align}\label{13.6}
	\frac{d}{dt} \io \frac{(\ueps+1)^p}{\del-\veps}
	&= -p(p-1)\io \f{(\ue+1)^{p-2}|\na \ue|^2}{\del-\ve} - 2p\io \f{(\ue+1)^{p-1}}{(\del-\ve)^2}\na \ue\cdot\na \ve \nn\\
	&+p(p-1)\io \f{(\ue+1)^{p-2}\ue}{(\del-\ve)(1+\eps\ue)^2}\na \ue\cdot\na \ve 
	+ p\io \f{(\ue+1)^{p-1}\ue}{(\del-\ve)^2(1+\eps\ue)^2}|\na \ve|^2 \nn\\
	&-2\io \f{(\ue+1)^p}{(\del-\ve)^3}|\na \ve|^2 - \io \f{(\ue+1)^p \ue\ve}{(\del-\ve)^2(1+\eps\ue)}.
  \end{align}
  Here, the elementary inequalities $1\le \f{\del}{\del-\ve}$, $(\ue+1)^{p-2}\ue \le (\ue+1)^{p-1}$ and 
  $\f1{(1+\eps\ue)^2}\le 1$ together with nonpositivity of the last term and an application of Young's inequality 
  due to \eqref{13.4} result in
  \begin{align*}
   	\frac{d}{dt} \io \frac{(\ueps+1)^p}{\del-\veps}
	&\le -p(p-1)\io \f{(\ue+1)^{p-2}|\na \ue|^2}{\del-\ve} 
	+ (2p+\del p(p-1))\io \f{(\ue+1)^{p-1}}{(\del-\ve)^2} |\na \ue||\na \ve|\\
	&+ (p\del-2)\io \f{(\ue+1)^p}{(\del-\ve)^2}|\na \ve|^2\\
   	&\le -\f{p(p-1)}4\io \f{(\ue+1)^{p-2}|\na \ue|^2}{\del-\ve} 
	+\kl{\f{(2p+\del p(p-1))^2}{3p(p-1)}+p\del-2} \io \f{(\ue+1)^p}{(\del-\ve)^2}|\na \ve |^2\\
   	&\le -\f{p(p-1)}{4\del} \io (\ue+1)^{p-2}|\na \ue|^2\quad \mbox{in } \Om\times (T^{(1)},\infty)
	\qquad \mbox{for all } \eps\in (\eps_j)_{j\in\N} \cap (0,\eps^{(1)}).
  \end{align*}
  A variant of the Gagliardo--Nirenberg inequality (cf. eg., \cite[Lemma 2.3]{li_lankeit_hapto}), stating that with suitable $c_1>0$
  we have
  \[
   	\norm[\Lom2]{\psi}\le c_1 \norm[\Lom2]{\na \psi} ^{\theta}\norm[\Lom{\f2p}]{\psi}^{1-\theta} 
	+ c_1 \norm[\Lom{\f2p}]{\psi} \qquad \forall \psi\in \Wom12\cap \Lom{\f2p}, 
  \]
  where $\theta=\f{n(p-1)}{n(p-1)+2}$, proves the existence of $c_2>0$ such that every $\psi\in \Wom12$ obeying an estimate of the 
  form $\norm[\Lom{\f2p}]{\psi}\le B$ with $B>0$ satisfies 
  \[
   \norm[\Lom2]{\psi}^{\f2{\theta}}\le c_2 B^{\f{2(1-\theta)}{\theta}} \norm[\Lom2]{\na \psi}^2 + c_2B^{\f2{\theta}}.
  \]
  An application to $\psi=(\ue+1)^{\f p2}$ on account of the consequence 
  \[
  	\io \f{(\ue+1)^p}{\del-v}\le \f2{\del}\io \psi^2 
  \]
  of $\ve \le \f{\del}2$ and \eqref{mass} shows that if for 
  $\eps \in (\eps_j)_{j\in\N} \cap (0,\eps^{(1)})$ and $t>T^{(1)}$ we let
  $\yeps(t):=\io \frac{(\ueps(\cdot,t)+1)^p}{\del-\veps(\cdot,t)}$
  and $\heps(t):=\io (\ueps(\cdot,t)+1)^{p-2} |\na\ueps(\cdot,t)|^2$,
  then with some positive constants $c_3,c_4$ and $c_5$,
  \be{13.7a}
  	\yeps'(t) + c_3 \heps(t) \le 0
  \ee
  and
  \be{13.7}
	\yeps'(t) + c_4 \yeps^\lam(t) \le c_5
  \ee
  for all $t>T^{(1)}$ and $\eps\in (\eps_j)_{j\in\N} \cap (0,\eps^{(1)})$,
  where $\lambda=\f{1}{\theta}=\f{n(p-1)+2}{n(p-1)}>1$. From \eqref{13.7}, we obtain \eqref{13.1} with 
  $T_{\star\star\star}(p)=T^{(1)}+1$ and $\eps_{\star\star\star}(p)=\eps^{(1)}$ by Lemma~\ref{lem:ode}. 
  Afterwards, \eqref{13.2} follows from an integration of \eqref{13.7a} with respect to time.
\qed 
As a first result of these estimates, we note that classical regularity theory turns them into even higher regularity information:
\begin{lem}\label{lem14}
  There exist $T_0>0, \eps_0\in (0,1)$ and $C>0$ such that
  \be{14.1}
	\|\ueps\|_{C^{2+\theta,1+\frac{\theta}{2}}(\bom\times [t,t+1])} \le C
	\qquad \mbox{for all $t>T_0$ and } \eps\in (\eps_j)_{j\in\N} \cap (0,\eps_0),
  \ee
  and that
  \be{14.2}
	\|\veps\|_{C^{2+\theta,1+\frac{\theta}{2}}(\bom\times [t,t+1])} \le C
	\qquad \mbox{for all $t>T_0$ and } \eps\in (\eps_j)_{j\in\N} \cap (0,\eps_0).
  \ee
\end{lem}
\newcommand{\holdernorm}[2][]{\|#2\|_{C^{#1,\f{#1}2}(\Ombar\times[t,t+1])}}
\proof
  With sufficiently large $p>3$ and the outcome of Lemma~\ref{lem13}, the H\"older regularity results 
  of \cite[Thm. 1.3]{porzio_vespri} provide $\eps$-independent bounds 
  for $\holdernorm[\theta_1]{\ve}$ and, subsequently, of $\holdernorm[\theta_2]{\ue}$ 
  for $t\in(T_{\star\star\star}(p)+1,\infty)$ with some $\theta_i\in (0,1)$, $i\in\{1,2\}$. 
  By \cite[Thm. IV.5.3]{LSU}, \eqref{14.2} follows, whenceupon \eqref{14.1} results from successive applications 
  of \cite[Thm.~1.1]{lieberman_gradient}, which combined with the previous estimates yield $\theta_3\in (0,1)$ and an 
  $\eps$-independent bound for $\holdernorm[\theta_3]{\del\ve\f{\ue}{(1+\eps\ue)^2}+\na \ve\cdot\na (\f{\ue}{(1+\eps\ue)^2}}$, 
  and \cite[Thm.~IV.5.3]{LSU}.
\qed
Both \eqref{13.2} and Lemma~\ref{lem14} have consequences for the solution $(u,v)$ constructed in Lemma~\ref{lem9}.
\begin{cor}\label{cor98}
  There are $T>0$ and $C>0$ such that for every $t\ge T$
  \begin{equation}\label{holdernorm-limit}
  	\holdernorm[2+\theta]{u}\le C
	\qquad \mbox{and} \qquad
	\holdernorm[2+\theta]{v}\le C;
  \end{equation}
  in particular, $(u,v)\in (C^{2,1}(\Ombar\times[T,\infty)))^2$, and that, moreover, 
  \begin{equation}\label{nablau2}
  	\int_T^\infty \io |\na u|^2 \le C.
  \end{equation}
\end{cor}
\proof
 The Hölder bounds of \eqref{14.1} show that for each $t\ge T=T_0$, $\ue\to u$ in $C^{2,1}(\Ombar\times[t,t+1])$ as $\eps=\eps_{j_k}\searrow 0$ along a suitable subsequence of $(\eps_j)_{j\in \N}$. Already pointwise convergence respects upper bounds on Hölder norms, therefore the first part of \eqref{holdernorm-limit} directly follows from \eqref{14.1} (the second from \eqref{14.2}, exactly the same way). 

 As to \eqref{nablau2}, we take any $p>3$ and observe that according to \eqref{13.2}, 
  $\{\na \ue^{\f p2} \mid \eps\in (\eps_j)_{j\in\N} \cap (0,\eps_{\star\star\star}(p))\}$ is bounded in $X:= L^2(\Om\times(\max\{T_{\star\star\star}(p),T_0\},\infty))$ and hence $\na \ue^{\f p2}\wto \na u^{\f p2}$ in $X$ (where identification of the limit relies on the $C^{2,1}$ convergence just derived). Weak sequential lower semicontinuity of the $L^2$ norm shows that hence
 \[
  	\int_T^\infty \io |\na u|^2\le \int_T^\infty \io (u+1)^{p-2}|\na u|^2 \le \liminf_{j\to \infty}  \int_T^\infty \io (u_{\eps_j}+1)^{p-2}|\na u_{\eps_j}|^2 <\infty.
 \]
\qed 
In deriving Theorem \ref{theo15} from this, 
we shall make use of the following observation made in \cite[Lemma 3.1]{taowin32}.
\begin{lem}\label{lem99}
  Let $t_0\in\R$ and $w:\bom\times [t_0,\infty) \to\R$ be uniformly continuous and such that
  \bas
	\int_t^{t+1} \io |w|^q \to 0
	\qquad \mbox{as } t\to\infty
  \eas
  with some $q>0$. Then
  \bas
	w(\cdot,t)\to 0
	\quad \mbox{in } L^\infty(\Om)
	\qquad \mbox{as } t\to\infty.
  \eas
\end{lem}
\proofc of Theorem \ref{theo15}. \quad
  With $T>0$ taken from Corollary \ref{cor98}, from \eqref{holdernorm-limit} we directly obtain (\ref{15.1}),
  and that $u$ is uniformly continuous in $\bom\times [T,\infty)$
  Since a Poincar\'e inequality yields $c_1>0$ such that
  \[
 	\int_T^\infty \norm[\Lom2]{u(\cdot,t)-\overline{u_0}}^2 dt \le c_1 \int_T^\infty \io |\na u|^2,
  \]
  in view of (\ref{nablau2}) we may draw on Lemma \ref{lem99} to infer that $u(\cdot,t)-\overline{u_0} \to 0$
  in $L^\infty(\Om)$ as $t\to\infty$.
  As, moreover,  Lemma \ref{lem12} in combination with \eqref{9.3} 
  implies that $\|v(\cdot,t)\|_{L^\infty(\Om)} \to 0$ as $t\to\infty$, the proof is complete.
\qed

\section{Acknowledgements}
The second author acknowledges support of the {\em Deutsche Forschungsgemeinschaft} (Project No.~462888149).

{\footnotesize 
  \setlength{\parskip}{0pt}
  \setlength{\itemsep}{0pt plus 0.05ex}
\setlength{\baselineskip}{0.92\baselineskip}
 \bibliographystyle{abbrv}

 }

\end{document}